\documentclass[%
english,%
twoside
]{%
amsart%
}

\usepackage{amsmath}

\usepackage{amssymb}

\usepackage{amstext}

\newtheorem{statement}{}
\newcommand{\bs}{\begin{statement}}
\newcommand{\es}{\end{statement}}

\newtheorem{example}{Example}

\newcommand{\bp}{\begin{proof}}

\newcommand{\ep}{\end{proof}}

\newcommand{\be}{\begin{equation}}
\newcommand{\ee}{\end{equation}}
\newcommand{\bea}{\begin{eqnarray}}
\newcommand{\eea}{\end{eqnarray}}
\newcommand{\beann}{\begin{eqnarray*}}
\newcommand{\eeann}{\end{eqnarray*}}
\newcommand{\bv}{\begin{verbatim}}
\newcommand{\ev}{\end{verbatim}}

\newcommand{\Z}{{\mathbb Z}}

\newcommand{\opon}{\ltimes}

\DeclareMathOperator{\abs}{abs}

\DeclareMathOperator{\Ad}{Ad}

\DeclareMathOperator{\Aut}{Aut}

\DeclareMathOperator{\card}{card}

\DeclareMathOperator{\GL}{GL}

\DeclareMathOperator{\Hom}{Hom}

\DeclareMathOperator{\Lin}{Lin}

\DeclareMathOperator{\slin}{sl}

\begin{document}

\title{%
	The Alexander- and Jones-invariants
	and the Burau module
      }

\author{Florin Constantinescu}

\address{%
	Institut f\"ur Angewandte Mathematik
\\
	Johann Wolfgang Goethe Universit\"at
\\
	Robert Mayer Str.~6-10
\\
	D-60054 Frankfurt am Main
\\
	Germany
\\
	email:	'constantinescu@ma\-the\-ma\-tik.uni-frankfurt.d400.de'
}

\author{Mirko L\"udde}

\address{%
	Institut f\"ur Reine Mathematik
\\
	Humboldt Universit\"at zu Berlin
\\
	Ziegelstrasse 13a
\\
	D-10099 Berlin
\\
	Germany
\\
	email:'luedde@mathematik.hu-berlin.de'
	}

\date{\today }

\begin{abstract}
	From the braid-valued Burau module over the braid group
	we construct the Yang-Baxter matrices
	yielding the Alexander- and the Jones knot invariants.
	This generalises an observation of V. F. R. Jones.
\end{abstract}

\maketitle

\section{%
	Introduction}
It has been known for long that the Burau representation
of Artin's braid group $ B_n $
can be used to construct
the Alexander polynomial invariant of knots.
There are at least two ways to accomplish this.
The topological approach uses the fact
that the first relative homology of the
cyclic covering of the knot's complement
has a presentation as a
$ \Z [t, t^-] $ module determined by the Burau matrices,
cf.~\cite{BurdeZieschang}.
Another approach,
the one which is generalised in this article,
constructs a solution of the Yang-Baxter equation
starting with the Burau representation.
This proceeds by extending the
$ 3 $-dimensional Burau representation
of $ B_3 $ to the
$ 2^3 $-dimensional Grassman algebra
of the representation space.
The Yang-Baxter matrix can then be turned into the
Alexander invariant,
e.~g.~by a state model on knot diagrams,
cf.~\cite{Jones1991,Kauffman}.
This idea goes back to Jones and has been investigated in
\cite{Kauffman,KauffmanSaleur1992}.

In this article,
following a proposal of \cite{Constantinescu19953},
we will obtain
the simplest Yang-Baxter solution associated to the deformation
$ U_t(\slin (2)), $
and therefore the Jones invariant.
We use a generalisation
of the Burau representation
that we called the
``braid-valued Burau representation''
in \cite{ConstantinescuLuedde1992,Luedde1992}.
By this term we mean the module obtained
as the relative augmentation ideal of the
free group $ F_n $ of rank $ n $
in the integral group ring
$ \Z [ B_n \opon F_n ]. $
Tensor products of these modules
can be suitably reduced:
an antisymmetrisation of tensor products
yields the Grassman algebra carrying the
``classical'' Burau representation
and therefore the Alexander invariant.
By similar but different relations
the Yang-Baxter matrix for the
Jones invariant is obtained.

Acknowledgements:
M.L.~thanks the
{\sl Graduiertenkolleg ``Geometrie und Nichtlineare Analysis''}
of the Humboldt University at Berlin for support.

\section{%
	The Burau representation
}
The definitions and facts on the braid group
used here
can be found e.g.~in \cite{Birman,BurdeZieschang},
if not stated otherwise.
The braid group
$ B_{n} $ on $ n $ strings
(over the Euclidean plane)
is the
group generated by the set
$
\{ \tau _{i}; i \in \{ 1, \ldots , n-1 \} \}
$
according to the relations of Artin,
$
\tau _{i} \tau _{j} =
\tau _{j} \tau _{i},
$ if
$ \abs (i-j) \geq 2, $
$
\tau _{i} \tau _{1+i} \tau _{i} =
\tau _{1+i} \tau _{i} \tau _{1+i}.
$

The braid group $ B_n $
faithfully acts onto the free group
$ F_n := < f_1, \ldots , f_n > $
of rank $ n. $
There is a monomorphism
$
\psi \in \Hom (B_n, \Aut (F_n))
$
(where we let the automorphisms act from the right)
defined by
$$
\psi(\tau _{i}) : f_{j} \mapsto
\left\{
       \begin{array}{cl}
       f_{i} f_{i+1} f_{i}^{-1},	& j = i,
\\
       f_{i},				& j = i + 1,
\\
       f_{j},				& j \not \in \{ i, i+1 \}
       \end{array}
\right. .
$$
We will only need the fact that $ \psi $ is a (anti-)homomorphism,
which can be checked by computation.

Now we can define the semidirect product
$ B_n F_n := B_n \opon _{\psi } F_n $
as the set $ B_n \times F_n $ with
multiplication
$ ( \alpha , f) ( \beta , g ) := ( \alpha \beta, (\psi(\beta))(f) g), $
which we will write simply as
$ \alpha f \beta g = \alpha \beta \beta (f) g. $

There are several approaches to the classical Burau representation.
Following W.~Magnus, cf.~\cite{Jones1991,Magnus1974},
one may investigate the $ B_n $ action on
abelianised groups $ U / [U,U] $ for suitable $ U \leq F_n. $
A different method proceeds by using the Fox-derivative
on the free group, cf.~\cite{Birman}.
A topological construction
via a flat connexion on a homology bundle
is obtained as a particular case in
\cite{Atiyah1989,Lawrence1990}.
The approach we are using
presumably had not been noted in the literature until
\cite{ConstantinescuLuedde1992,Luedde1992}.
For an algebraic derivation of all these methods
and their generalisations see \cite{Luedde19952}.
\begin{statement}[Braid-valued Burau module]
The map
$$
\tau _i \mapsto
\tau _i
\left(
\begin{array}{ccccc}
{\bf 1}_{i-1} &0	&0	&0
\\
0	&(1 - f_{i} f_{i+1} f_{i}^{-1})	&f_{i}	&0
\\
0	&1		&0			&0
\\
0	&0		&0		&{\bf 1}_{n-i-1}
\end{array}
\right)
$$
uniquely extends to a monomorphism
$ B_n \rightarrow \GL (n,\Z [ B_n F_n ]). $
\end{statement}
\bp
The relative augmentation ideal
$
\Lin _{\Z B_nF_n} \{ (f_i - 1); i \in \{ 1, \ldots , n \}  \}
$
of the free group in the integral group ring
$ \Z [ B_n F_n ] $
is free of rank $ n $
as a left $ B_n F_n $ module over the set
$ \{ s_i := (f_i - 1); i \in \{ 1, \ldots , n \} \}. $
As an ideal, by multiplication from the right,
it is a module over $ B_n, $
$
(f_j - 1) \mapsto
(f_j - 1) \tau _i = \tau _i (\tau _i(f_j) - 1).
$
The element
$ \tau _i(f_j) := \psi (\tau _i)(f_j) $
is determined by Artin's action,
so we obtain the equation
$
s_j \tau _i =
\tau _{i}
\left\{ \begin{array}{ll}
(1 - f_{i}f_{i+1}f_{i}^{-1})s_{i} + f_{i}s_{i+1}, & j=i
\\
s_{i},						& j = i + 1
\\
s_{j},						& j \not \in \{ i, i+1 \}
\end{array} \right. .
$
This action is faithful, since the matrix representative of
$ \alpha \in B_n $ has the form
$ \alpha S, $
with $ S $ an $ n $ by $ n $ matrix over the ring
$ \Z F_n \hookrightarrow \Z B_n F_n. $
\ep

\begin{example}[Burau module]
The ring homomorphism
$ \Z [B_n F_n] \rightarrow \Z [t, t^-], $
with an indeterminate $t$,
defined by
$ \tau _i \mapsto 1, $
$ f_j \mapsto t, $
applied to the braid-valued Burau matrices,
yields the representation
$$
\tau _i \mapsto
\left(
\begin{array}{ccccc}
{\bf 1}_{i-1}&0		&0	&0
\\
0	& 1 - t 	&t	&0
\\
0	&1		&0	&0
\\
0	&0		&0	&{\bf 1}_{n-i-1}
\end{array}
\right) .
$$
\end{example}

\section{%
	Construction of Yang-Baxter matrices
}
Having at hand the Burau representation,
we recall Jones' construction of a
Yang-Baxter matrix from it.
This prepares us for the general procedure.
\bs[Alexander invariant from Burau module]
Let
$ R := \Z [t,t^-] $
be the ring of
Laurent polynomials in an indeterminate $ t. $
Let
$ \rho \in \Hom (B_3, \Aut (R^3)) $
be a Burau representation,
given by the matrices
$
\rho _1 := \left( \begin{array}{cc}
B & 0 \\
0 & 1
\end{array} \right)
$
and
$
\rho _2 := \left( \begin{array}{cc}
1 & 0 \\
0 & B
\end{array} \right) ,
$
where
$
B := \left( \begin{array}{cc}
1 - t & t \\
1 & 0
\end{array} \right) .
$
Then the natural extension of $ \rho $
to the exterior algebra
$ \Lambda (R^3) $
(where the $ \rho _i $ act as algebra homomorphisms)
is isomorphic to a Yang-Baxter representation
$ \Upsilon \in \Aut (V \otimes V) $
with a $ 2 $-dimensional free $ R $ module $ V. $
\es
\bp
The statement and proof are taken from
\cite{Kauffman}, sect.~I.13, pp.~208.
Let $ R^3 $ as an $R$ module have the basis
$ v_1 = (1,0,0), $
$ v_2 = (0,1,0), $
$ v_3 = (0,0,1). $
The representation $ \rho $ acts by multiplication from the right
onto these row vectors.
Let $ V $ be the free left $ R $ module with basis
$ e_1, $ $ e_2. $
Then we define an isomorphism of free left $R$ modules
$ \phi \in \Hom (\Lambda (R^3), V \otimes V \otimes V) $
by sending
$
( 1, v_1, v_2, v_1 \wedge v_2,
v_3, v_1 \wedge v_3, v_2 \wedge v_3, v_1 \wedge v_2 \wedge v_3 )
$
to
$
(
e_1 e_1 e_1, e_1 e_2 e_1, e_2 e_1 e_1, e_2 e_2 e_1,
e_1 e_1 e_2, e_1 e_2 e_2, e_2 e_1 e_2, e_2 e_2 e_2
).
$
Computing
$ \sigma _i = \phi \circ \rho _i \circ \phi ^{-1} $
we find,
$ \sigma _1 = \Upsilon \otimes 1 $
and
$ \sigma _2 = 1 \otimes \Upsilon $
with the matrix
$$
\Upsilon :=
\left( \begin{array}{cccc}
1 & 0 & 0 & 0
\\
0 & 1- t & t & 0
\\
0 & 1 & 0 & 0
\\
0 & 0 & 0 & -t
\end{array} \right)
$$
(where the rows correspond to the coefficients of the images
$ \Upsilon (e_i \otimes e_j) $
in the ordered basis
$ e_1e_1, $ $ e_1e_2, $ $ e_2e_1, $ $ e_2e_2. $)
This matrix satisfies the (permuting form of the)
Yang-Baxter equation,
$
\sigma _1 \sigma _2 \sigma _1 = \sigma _2 \sigma _1 \sigma _2,
$
as a consequence of the braid relations obeyed by
$ \rho _1 $ and $ \rho _2. $
\ep

A rescaling and a
change of basis transforms the matrix
into the one
that via a state model on knot diagrams
yields the Alexander invariant,
as described in \cite{Kauffman}, sect.~I.12, pp.~174.

We want to apply a similar technique to the
braid-valued Burau module.
By Artin's combed normal form for the braid group,
the group $ B_n F_n $ can be imbedded into
$ B_{1+n}. $
These imbeddings can be iterated to build the groups
$ B_{n,j} := B_n \opon F^{(1)} \opon \cdots \opon F^{(j)}, $
where $ F^{(j)} := F_{n+j-1}. $
We need to know two facts on this imbedding.
The first is that the generators of $ B_n $ are mapped as
$ \tau _i^{(n)} \mapsto \tau _i^{(1+n)} $
for $ i \in \{ 1, \ldots , n-1 \} , $
where superscripts indicate the respective groups.
The second is, in which way the image of the free group
$ F^{(j)} := <f_1^{(j)}, \ldots , f_{n+j-1}^{(j)}> $
in $ B_{n+j} \hookrightarrow \Aut (F^{(l)}) $
acts onto the generators of $ F^{(l)} $
for $ l > j. $
This action is given by
($ \epsilon \in \{ -1, 1 \} $)
$$
f_{i}^{(j) \epsilon }(f_k^{(l)}) =
\left\{ \begin{array}{ll}
f_k^{(l)},
	& k < i \text{ or } n+j-1 < k
\\
\Ad ((f_i^{(l)}f_{n+j-1}^{(l)})^{\epsilon })(f_k^{(l)}),
	& k \in \{ i, n+j-1 \}
\\
\Ad ([f_{n+j-1}^{(l) \epsilon },f_i^{(l) \epsilon }]^{-\epsilon })
(f_k^{(l)}),
	& i < k < n+j-1
\end{array} \right. ,
$$
where
$ k < n+l-1, $
$ i < n+j-1, $
$ \Ad (x)(y) := xyx^-, $
$ [x, y] := xyx^-y^-. $
These equations are equivalent to the relations for the
generators of the pure braid group.
We consider the groups
$ B_{3,j} $
for $ j \in \{ 1, 2, 3 \} . $
Let
$
I^{(j)} := \Lin _{\Z B_{3,j}} \{ s_i^{(j)} := (f_i^{(j)} - 1) \}
$
be the relative augmentation ideal of $ F^{(j)} $
in the ring $ \Z B_{3,j}.$
Furthermore,
let a left $B_{3,3}$ right $B_3$ bimodule be defined
as the sum of tensor products
$$
M :=
\Z [B_{3,3}] {\mathbf 1} \oplus
I^{(3)} \oplus
I^{(3)} \otimes _{B_{3,2}} I^{(2)} \oplus
I^{(3)} \otimes _{B_{3,2}} I^{(2)} \otimes _{B_{3,1}} I^{(1)}.
$$
Define a right $B_{3,3}$ module structure on
the ring $ \Z [t,t^-] $ of Laurent polynomials
by mapping the generators of $ B_{3,3} $ as
$ \tau _i \mapsto 1, $
for $ i \in \{ 1, 2 \} $
and
$ f_k^{(j)} \mapsto t. $
Guided by Jones' construction of the representation on
the Grassman algebra,
we construct quotients of rank $ 2^3 $ of $ M. $
\bs[Invariants from braid-valued Burau module]
The representation of $B_3$ defined by the
Yang-Baxter matrix
$$
R :=
\left( \begin{array}{cccc}
1 & 0 & 0 & 0
\\
0 & 1 - t & t & 0
\\
0 & 1 & 0 & 0
\\
0 & 0 & 0 & 1
\end{array} \right) ,
$$
as well as the representation by the previously defined
matrix $ \Upsilon $
can be obtained as suitable quotients from the tensor product
$ \Z [t,t^-] \otimes _{B_{3,3}} M $
regarded as left $ \Z [t,t^-] $ right $ B_3 $ bimodule.
\es
\bp
We compute the right action of $ \tau _1 $ induced
by Artin's automorphisms on some
particular basis elements of $ M. $
For notational convenience we drop the tensor product
symbol and set
$ a_i^{(j)} := ( 1 - f_i^{(j)} f_{1+i}^{(j)} f_i^{(j)-} ). $
We get
$
{\mathbf 1} \tau _1 =
\tau _1 {\mathbf 1},
$
$
s_1^{(3)} \tau _1 =
\tau _1 ( a_1^{(3)} s_1^{(3)} + f_1^{(3)} s_2^{(3)} ),
$
$
s_2^{(3)} \tau _1 =
\tau _1 s_1^{(3)},
$
$
s_3^{(3)} \tau _1 =
\tau _1 s_3^{(3)},
$
$
s_1^{(3)} s_2^{(2)} \tau _1 =
\tau _1 ( a_1^{(3)} s_1^{(3)} + f_1^{(3)} s_2^{(3)} ) s_1^{(2)},
$
$
s_1^{(3)} s_3^{(2)} \tau _1 =
\tau _1 ( a_1^{(3)} s_1^{(3)} + f_1^{(3)} s_2^{(3)} ) s_3^{(2)},
$
$
s_2^{(3)} s_3^{(2)} \tau _1 =
\tau _1 s_1^{(3)} s_3^{(2)}
$
and finally,
$
s_1^{(3)} s_2^{(2)} s_3^{(1)} \tau _1 =
\tau _1 ( a_1^{(3)} s_1^{(3)} + f_1^{(3)} s_2^{(3)} )
s_1^{(2)} s_3^{(1)}.
$
Similarly, the action of the second generator $ \tau _2 $
is as follows:
$
{\mathbf 1} \tau _2 =
\tau _2 {\mathbf 1},
$
$
s_1^{(3)} \tau _2 =
\tau _2 s_1^{(3)},
$
$
s_2^{(3)} \tau _2 =
\tau _2 ( a_2^{(3)} s_2^{(3)} + f_2^{(3)} s_3^{(3)} ),
$
$
s_3^{(3)} \tau _2 =
\tau _2 s_2^{(3)},
$
$
s_1^{(3)} s_2^{(2)} \tau _2 =
\tau _2
( a_2^{(2)} s_1^{(3)} s_2^{(2)} + f_2^{(2)} s_1^{(3)} s_3^{(2)}),
$
$
s_1^{(3)} s_3^{(2)} \tau _2 =
\tau _2 s_1^{(3)} s_2^{(2)},
$
$
s_2^{(3)} s_3^{(2)} \tau _2 =
\tau _2
( a_2^{(3)} s_2^{(3)} + f_2^{(3)} s_3^{(3)} ) s_2^{(2)},
$
$
s_1^{(3)} s_2^{(2)} s_3^{(1)} \tau _2 =
\tau _2 ( a_2^{(2)} s_1^{(3)} s_2^{(2)} s_2^{(1)}
+ f_2^{(2)} s_1^{(3)} s_3^{(2)} s_2^{(1)}).
$
We now pass to the quotient of
$ \Z [t,t^-] \otimes _{B_{3,3}} M $
as a left $ \Z [t, t^-] $ module
by imposing the relations
(note the similarity with the relations of a
Grassman algebra)
\begin{eqnarray*}
s_i^{(3)} s_j^{(2)} & = &
	\left\{ \begin{array}{rl}
	0,
	& i = j
	\\
	t^{-} s_j^{(3)} s_i^{(2)},
	& i > j
	\end{array}\right. ,
\\
s_{i}^{(3)} s_j^{(2)} s_k^{(1)} & = &
	\left\{ \begin{array}{rl}
	0, & \card \{ i, j, k \} < 3
	\\
	t^{-} s_{j}^{(3)} s_{i}^{(2)} s_{k}^{(1)},
	& i > j
	\\
	t^{-} s_{i}^{(3)} s_{k}^{(2)} s_{j}^{(1)},
	& j > k
	\end{array} \right. .
\end{eqnarray*}
This quotient $ Q $ is
a free left $ \Z [t,t^-] $ module
with basis given by the $ 2^3 $ elements
$
( 1, s_1, s_2, s_1 s_2, s_3, s_1 s_3, s_2 s_3, s_1 s_2 s_3 ) .
$
We obtain an induced action of
$ \tau _1 $ and $ \tau _2 $
on $ Q, $
$
{\mathbf 1} \tau _1 = {\mathbf 1},
$
$
s_1^{(3)} \tau _1 = (1 - t) s_1^{(3)} + t s_2^{(3)},
$
$
s_2^{(3)} \tau _1 = s_1^{(3)},
$
$
s_3^{(3)} \tau _1 = s_3^{(3)},
$
$
s_1^{(3)} s_2^{(2)} \tau _1 = s_1^{(3)} s_2^{(2)},
$
$
s_1^{(3)} s_3^{(2)} \tau _1 =
( ( 1 - t ) s_1^{(3)} + t s_2^{(3)} ) s_3^{(2)},
$
$
s_2^{(3)} s_3^{(2)} \tau _1 = s_1^{(3)} s_3^{(2)},
$
$
s_1^{(3)} s_2^{(2)} s_3^{(1)} \tau _1 =
s_1^{(3)} s_2^{(2)} s_3^{(1)}
$
and similarly for the action of $ \tau _2. $
In order to obtain a Yang-Baxter representation of $ B_3, $
consider the free left $ \Z [t,t^-] $ module
$ V \otimes V \otimes V $
of rank $ 2^3 $ with
$ V := \Z [t, t^-]^{2} $
and define an isomorphism
$ Q \rightarrow V \otimes V \otimes V $
by sending the ordered basis above to
$
( e_1 e_1 e_1, e_1 e_2 e_1, e_2 e_1 e_1, e_2 e_2 e_1,
 e_1 e_1 e_2, e_1 e_2 e_2, e_2 e_1 e_2, e_2 e_2 e_2 ).
$
On $ V \otimes V \otimes V $
the braid generators $ \tau _1 $ and $ \tau _2 $
are then found to act by matrices
$ R \otimes 1 $ and $ 1 \otimes R, $
respectively.

Finally we notice that a similar construction,
where we impose the exact Grassman relations,
leads to the matrix $ \Upsilon $ precisely as in
the previous lemma.
\ep

By a rescaling and a change of basis in $ V \otimes V, $
the Yang-Baxter matrix $ R $ is found to be
the universal R-matrix of $U_t(\slin (2))$ in its
fundamental representation,
cf.~\cite{ChariPressley}, ex.~6.4.12, pp.~205.
So either by a state model,
see \cite{Kauffman}, sect.~I.11, pp.~161,
or by Turaev's theorem,
see \cite{ChariPressley}, sect.~15.2, pp.~504,
the Jones polynomial can be obtained.

\bibliographystyle{amsalpha}

\providecommand{\bysame}{\leavevmode\hbox to3em{\hrulefill}\thinspace}

\end{document}